\documentclass{tac}
\usepackage{amssymb}
\usepackage{bbold}
\input diagxy

\def\lax{{\rm Lax}_{\rm N}}
\def\cat{\mathbb{C}{\rm at}}
\def\pos{\mathbb{P}{\rm os}}
\def\loc{\mathbb{L}{\rm oc}}
\def\spaces{\mathbb{T}{\rm op}}

\def\rel{\mathbb{R}{\rm el}}
\def\tov{\hbox{$\to<150>$\hskip -.2in \raisebox{1.7pt}{\tiny$\bullet$} \hskip .1in}}
\def\id{{\rm id}}
\def\idv{\id^{\hbox{\tiny$\bullet$}}}

\def\exp{exponentiable}
\def\expty{exponentiability}

\title{Exponentiability via Double Categories}
\author{Susan Niefield}
\address{Union College \\ Department of Mathematics \\ Schenectady, NY 12308}
\eaddress{niefiels@union.edu}

\begin{document} 
\maketitle

\begin{abstract} 
For a small category  $B$ and a double category $\mathbb D$, let ${\rm Lax}_N(B,\mathbb D)$ denote the category whose objects are vertical normal lax functors $B\to<125>\mathbb D$ and morphisms  are horizontal lax transformations.  It is well known that ${\rm Lax}_N(B,\mathbb C{\rm at})\simeq {\rm Cat}/B$, where $\mathbb{C}{\rm at}$ is the double category of small categories, functors, and profunctors.  In \cite{glue}, we generalized this equivalence to certain double categories, in the case where $B$ is a finite poset.  In \cite{Street}, Street showed that $Y\to<125> B$ is exponentiable in $ {\rm Cat}/B$ if and only if the corresponding normal lax functor $B\to<125> \mathbb C{\rm at}$ is a pseudo-functor.   Using our generalized equivalence, we show that a morphism $Y\to<125> B$ is exponentiable in $ {\mathbb D}_0/B$ if and only if the corresponding normal lax functor $B\to<125>\mathbb D$ is a pseudo-functor {\it plus} an additional condition that holds for all  $X\to<125> B$ in $\rm Cat$.  Thus, we obtain a single theorem which yields characterizations of certain exponentiable morphisms of small categories, topological spaces, locales, and posets.
\end{abstract} 

\section{Introduction}
Suppose $\cal A$ is a category with finite limits.  An object $Y$ of $\cal A$ is called {\it \exp} if the functor $-\times Y\colon {\cal A}\to<125>{\cal A}$ has a right adjoint, denoted by $(\ )^Y$.  A morphism is called {\it \exp} if it is {\exp} in ${\cal A}/Y$.

\smallskip
Exponentiable morphisms in the category $ {\rm Cat}$ of small categories were characterized independently by Giraud \cite{Giraud} and  Conduch\'{e} \cite{Conduche} as those functors satisfying a factorization lifting property now known as the {\it Giraud-Conduch\'{e} condition}.  Exponentiable morphisms in the category $\rm Top$ of topological spaces were characterized by the author in \cite{thesis} (see also \cite{cart,posets,perfect}) as those satisfying a somewhat technical condition  (see Lemma~\ref{cart} below) which was used to show that the inclusion of a subspace of $B$ is {\exp}  if and only if it is locally closed, and also to establish the {\expty} of perfect maps as well as locally compact spaces over a locally Hausdorff base. 

\smallskip
The obstruction to {\expty} in each of these two categories is quite different.  In $ {\rm Cat}$, the Giraud-Conduch\'{e} condition is used to define composition of morphisms in the category that serves as the exponential, and the unit and counit follow.  Whereas in $\rm Top$, one can always define a candidate for the exponential for which the unit is continuous, but the extra condition is needed for the continuity of the counit.

\smallskip
There is a more recent characterization of {\expty} in $ {\rm Cat}$.  In a 1986 hand-written manuscript, referenced in his 2001 unpublished note \cite{Street}, Street  used the equivalence (attributed to B\'{e}nabou) between $ {\rm Cat}/B$ and a category $\lax(B,Prof)$ to show that a functor $Y\to<125>B$ is {\exp} if and only if the corresponding normal lax functor $B\to<125>Prof$ is a pseudo-functor.   Objects of $\lax(B,Prof)$ are normal lax functors from $B$ to the bicategory $Prof $ of small categories and profunctors, 
and  morphisms of $\lax(B,Prof)$ are functor-valued lax transformations.  Thus, B\'{e}nabou's equivalence can by viewed as taking place in the realm of double categories (in the sense of \cite{Ehres} or \cite{GP99}).  In particular, for a double category $\mathbb D$, we can consider the category $\lax(B,\mathbb D)$ whose objects are vertical normal lax functors and morphisms are horizontal lax transformations.  

\smallskip
In \cite{glue}, we established the equivalence between  $\lax(B,\mathbb D)$ and $\mathbb D_0/\Gamma_B1$, for certain double categories $\mathbb D$, in the case where $B$ is a finite poset and the constant functor $\mathbb D_0\to<125>\lax(B,\mathbb D)$ admits a left adjoint $\Gamma_B$.  When  $B=\mathbb2$, we know $\Gamma_B1=\mathbb2$ in $\rm Cat$ and $\rm Pos$.  It is the Sierpinski space $\mathbb2$ in $\rm Top$ and the Sierpinski locale $S$ in $\rm Loc$.  The poset $B=\mathbb2$ was also used in \cite{glue} to define open and closed inclusions in $\mathbb D_0$ and obtain a general construction of exponentials for locally closed inclusions over an {\it arbitrary} base, which we then applied to  $ {\rm Cat}$, $\rm Top$, $\rm Loc$, and $\rm Pos$. 

\smallskip
In this paper, we characterize the {\exp} objects $Y\colon B\to<125>\mathbb D$ of $\lax(B,\mathbb D)$ when $B$ is the 3-element linearly-ordered poset, $\lax(B,\mathbb D)$ has finite limits, and $\mathbb D$ has, what we call, a {\it zero object}.  Double categories $\mathbb D$ with these properties include $\cat$, $\pos$, $\spaces$, $\loc$, and $\rel$.  In particular, we show that $Y$ is {\exp} if and only if  $-\times Y$ preserves pseudo-functors, $Y_b$ is {\exp}, for all $b\in B$, and  $Y_b\tov Y_c$ is {\exp} as an object of $\mathbb D_1$, for all $b<c$ in $B$.  Using the equivalence established in \cite{glue}, we thus obtain a characterization of {\expty} in $\mathbb D_0/\Gamma_B1$, for a general poset $B$, which applies to $ {\rm Cat}$, ${\rm Pos}$, $\rm Top$,  and $\rm Loc$,  with the additional assumption that $B$ is finite in the latter two cases.  Note that every vertical morphism  is {\exp} in $\mathbb D_1$, when $\mathbb D$ is $\cat$ or $\pos$.

\smallskip
We proceed as follows.  In Section~2, we recall the definition of double category, and introduce zero objects as well as the double categories that will be considered throughout.  The definition of $\lax(B,\mathbb D)$ and the characterization of its {\exp} objects are presented in Section~3.  We conclude, in Section~4, by characterizing {\exp} objects of $\mathbb D_1$, for the two remaining cases, namely, $\spaces$ and $\loc$.

\smallskip
The author would like to thank the referee for carefully reading the preliminary versions of this paper and providing useful suggestions and corrections.

\section{Double Categories}
A {\it double category} $\mathbb D$ is a weak internal category 
$$\bfig 
\morphism(0,0)<450,0>[\mathbb D_1\times_{\mathbb D_0}\mathbb D_1`\mathbb D_1;c]
\morphism(450,0)|a|/@{>}@<5pt>/<350,0>[\mathbb D_1`\mathbb D_0;d_0]
\morphism(800,0)|m|<-350,0>[\mathbb D_0`\mathbb D_1;\Delta]
\morphism(450,0)|b|/@{>}@<-5pt>/<350,0>[\mathbb D_1`\mathbb D_0;d_1]
\efig$$ 
in $\rm CAT$.  It consists of objects (those of $\mathbb D_0$), two types of morphisms: horizontal (morphisms of $\mathbb D_0$) and vertical (objects of $\mathbb D_1$ with domain and codomain given by $d_0$ and $d_1$), and cells (morphisms of $\mathbb D_1$) of the form
$$\bfig
\square<350,350>[X_0`Y_0`X_1`Y_1;f_0`m`n`f_1]
\place(0,200)[\hbox{\tiny$\bullet$}]
\place(350,200)[\hbox{\tiny$\bullet$}]
\place(175,175)[\to<100>]
\efig\eqno{(1)}$$
Composition of morphisms and identities are defined horizontally in $\mathbb D_0$ and vertically using $c$ and $\Delta$, respectively. Cell composition is given horizontally in $\mathbb D_1$ and vertically via $c$.  Cells in which the horizontal morphisms are identities are called {\it special cells}.  

\medskip
There are five double categories of interest in this paper.

\begin{description}

\item[{\rm (E1)}] $\cat$ has small categories as objects, functors and profunctors as horizontal and vertical morphisms, respectively, and natural transformations $m\to<125>n(f_0-,f_1-)$ as cells of the form $(1)$.

\item[{\rm(E2)}]  $\spaces$ has topological spaces as objects and continuous maps as horizontal morphisms.  Vertical morphisms $m\colon X_0\tov X_1$ are finite intersection-preserving maps $m\colon {\cal O}(X_0)\to<125> {\cal O}(X_1)$ on the open set lattices, and there is a cell of the form $(1)$ if and only if $f_1^{-1}n\subseteq mf_0^{-1}$.

\item[{\rm(E3)}] $\loc$ has locales as objects, locale morphisms (in the sense of \cite{SS}) as horizontal morphisms, and
 finite meet-preserving maps as vertical morphisms. There is a cell of the form $(1)$ if and only if $f_1^*n\le mf_0^*$.

\item[{\rm(E4)}] $\pos$ has partially-ordered sets as objects and order-preserving maps as horizontal morphisms.  Vertical morphisms $m\colon X_0\tov X_1$ are order ideals (i.e., up-sets) $m\subseteq  X_0^{op}\times X_1$, and there is a cell of the form $(1)$ if and only if $(x_0,x_1)\in m \Rightarrow (f_0(x_0),f_1(x_1))\in n$.

\item[{\rm(E5)}] $\rel$ has sets as objects, functions and relations as horizontal and vertical morphisms, respectively, and there is a cell of the form $(1)$ if and only if $(x_0,x_1)\in m \Rightarrow (f_0(x_0),f_1(x_1))\in n$.

\end{description}

Our most general result, Theorem~\ref{lax}, will follow from  properties shared by the five double categories.  Although (E1)--(E5) are all framed bicategories (in the sense of \cite{Shulman}) and the first four have $\mathbb 2$-glueing (in the sense of \cite{glue}), these conditions will not be used until we apply the main theorem to obtain the {\expty} results in (E1)--(E5).

\smallskip
An object $0$ of $\mathbb D$ is called a {\it zero object} if it is horizontally initial, vertically both initial and terminal, and there exists a unique cell
$$\bfig
\square/->``->`->/<350,450>[X_0`Y_0`X_1`Y_1;f_0``n`f_1]
\place(0,160)[\hbox{\tiny$\bullet$}]
\place(0,370)[\hbox{\tiny$\bullet$}]
\morphism(0,450)<0,-225>[X_0`0;]
\morphism(0,225)<0,-225>[0`X_1;]
\place(350,250)[\hbox{\tiny$\bullet$}]
\place(175,225)[\to<100>]
\efig$$
for all $f_0, f_1, n$.  Note that the latter condition implies that $X_0\tov 0\tov X_1$ is an initial object in the category of vertical morphisms $X_0\tov  X_1$ and special cells. 

\smallskip
The double categories in the five examples each have an initial object which is a zero object.  In fact, if $\mathbb D$ is any framed bicategory, then any horizontal initial object which is vertically initial and terminal is easily seen to be a zero object.

\section{Exponentiability of Normal Lax Functors}
Suppose $B$ is a poset and $\mathbb D$ is a double category.  A {\it vertical normal lax  functor} $X\colon B\to<125>\mathbb D$ consists of an object $X_b$ of $\mathbb D$, for every $b\in B$, a vertical morphism $X_{bc}\colon X_b\tov X_c$, for every $b<c$, and a special cell $X_{cd} X_{bc}\to<125> X_{bd}$, called a {\it comparison cell}, for every $b<c<d$, satisfying the usual coherence conditions.  In particular, we are assuming our normal lax functors are {\it strict} normal in the sense that $X_{bb}$ is the vertical identity morphism on $X_b$, for all $b$.  A normal lax functor for which the comparison cells are all isomorphisms is called a {\it pseudo-functor}.   A {\it horizontal lax transformation} $f\colon X\to<125> Y\colon  B\to<125>\mathbb D$ consists of a horizontal morphism $f_b\colon X_b\to<125>Y_b$, for all $b\in B$, and a cell
$$\bfig
\square<350,350>[X_b`Y_b`X_c`Y_c;f_b`X_{bc}`Y_{bc}`f_c]
\place(0,200)[\hbox{\tiny$\bullet$}]
\place(350,200)[\hbox{\tiny$\bullet$}]
\place(175,175)[\to<100>]
\efig$$
for every  $b<c$, compatible with the comparison cells for $X$ and $Y$.  Vertical normal lax functors and horizontal transformations form a category which we denote by $\lax(B,\mathbb D)$.  Note that $\lax(\mathbb1,\mathbb D)\cong\mathbb D_0$ and $\lax(\mathbb2,\mathbb D)\cong\mathbb D_1$.

\smallskip
If $\mathbb D$ has a zero object, then pseudo-functors $X\colon \mathbb 3\to<125> \mathbb D$ can be described as follows.  Given {\it any} normal lax functor   $X\colon \mathbb 3\to<125> \mathbb D$, there is a commutative diagram in $\lax(\mathbb 3,\mathbb D)$ of the form
$$\bfig
\morphism(300,800)<0,-200>[\scriptstyle0`\scriptscriptstyle{X_1};]
\place(300,715)[\cdot] 
\morphism(300,600)<0,-200>[\scriptscriptstyle{X_1}`\scriptstyle0;]
\place(300,520)[\cdot] 
\morphism(0,500)<0,-200>[\scriptscriptstyle{X_0}`\scriptscriptstyle{X_1};]
\place(0,415)[\cdot] 
\morphism(0,300)<0,-200>[\scriptscriptstyle{X_1}`\scriptstyle0;]
\place(0,220)[\cdot] 
\morphism(600,500)<0,-200>[\scriptstyle0`\scriptscriptstyle{X_1};]
\place(600,415)[\cdot] 
\morphism(600,300)<0,-200>[\scriptscriptstyle{X_1}`\scriptscriptstyle{X_2};]
\place(600,220)[\cdot] 
\Atriangle(0,300)/->`->`/<300,300>[\scriptscriptstyle{X_1}`\scriptscriptstyle{X_1}`\scriptscriptstyle{X_1};``]
\Vtriangle(0,0)/`->`->/<300,300>[\scriptscriptstyle{X_1}`\scriptscriptstyle{X_1}`X;``]
\efig\eqno{(2)}$$
where the morphisms and cells are identities  or uniquely determined by the definition of zero object.  If $X$ is a functor, i.e., $X_{12}X_{01}=X_{02}$ and the comparison cell is the identity, then one can show that (2) is a pushout in $\lax(\mathbb 3,\mathbb D)$. Thus, we get:

\lemma\label{po} The diagram (2) is a pushout in $\lax(\mathbb 3,\mathbb D)$ if and only if $X$ is a pseudo-functor.
\endlemma

\proof Since (2) is a pushout when $X$ is a functor, it follows that such a diagram is a pushout if and only if there is an isomorphism from $X$ to the corresponding functor $X_0 \tov X_1\tov X_2$ such that the cells $X_{01}\to<125>X_{01}$ and $X_{12}\to<125> X_{12}$ are identities if and only if $X$ is a pseudo-functor.
\endproof

\lemma\label{subposet} If $\mathbb D$ has a zero object and $A$ is a subposet of $B$, then the restriction functor $(\ )_A\colon\lax(B,\mathbb D)\to<125> \lax(A,\mathbb D)$ has a left adjoint $L_A$ such that $(\ )_AL_A=id$.
\endlemma

\proof Given $X\colon A\to<125>\mathbb D$, define
$$(L_AX)_b=\cases{X_b& if $b\in A$\cr \ 0 & otherwise}$$
and given $b<c$, let $(L_AX)_{bc}=X_{bc}$, if $b,c\in A$, and let $(L_AX)_{bc}$ be the unique vertical morphism to or from $0$, otherwise.  That $L_AX$ is a normal lax functor follows directly from the definition of zero object, and the result easily follows.
\endproof

\smallskip
We will use the following notation for the functors $L_A$  and $(\ )_A$ in some special cases.  We write $L_b$ and $(\ )_b$ when $A=\{b\}$ and $L_{bc}$ when $A$ is the subposet with two elements $b<c$.  Similarly, we use the notation $L_{bcd}$ in the case where $b<c<d$.  We also write $L_b\colon\mathbb D_0\to<125>\mathbb D_1$  and $(\ )_b\colon\mathbb D_1\to<125>\mathbb D_0$ for the functors induced by the isomorphisms $\lax(\mathbb1,\mathbb D)\cong\mathbb D_0$ and $\lax(\mathbb2,\mathbb D)\cong\mathbb D_1$, where $b=0,1$.

\smallskip 
Suppose $\lax(B,\mathbb D)$ has finite products. Then so does  $\lax(A,\mathbb D)$, for all subposets $A\subseteq B$, since $(\ )_A$ has a left adjoint by Lemma~\ref {subposet}.  In particular, $\mathbb D_0$ has finite products and $(X\times Y)_b\cong X_b\times Y_b$ in $\mathbb D_0$, for all $b\in B$ and $X,Y\in \lax(B,\mathbb D)$.  

\smallskip
Note that if $\mathbb D_0$ has  chosen products, we do not necessarily know that we can take \break $(X\times Y)_b$ to be the chosen product.  However, we will assume $\mathbb D$ is {\it horizontally invariant}, in the sense of \cite{GP99}, and then this problem disappears.  The five double categories (E1)--(E5) of interest all have companions and conjoints, and hence are horizontally invariant.

\smallskip
Given $b<c<d$, we have cells
$$\bfig
\square/->``->`->/<550,550>[X_b\times Y_b`X_b\times Y_b`X_d\times Y_d`X_d\times Y_d;``X_{bd}\times Y_{bd}`]
\place(0,175)[\hbox{\tiny$\bullet$}]
\place(0,450)[\hbox{\tiny$\bullet$}]
\morphism(0,550)<0,-275>[X_b\times Y_b`X_c\times Y_c;X_{bc}\times Y_{bc}]
\morphism(0,275)<0,-275>[X_c\times Y_c`X_d\times Y_d;X_{cd}\times Y_{cd}]
\place(550,300)[\hbox{\tiny$\bullet$}]
\place(320,320)[\hbox{\small$\theta_{\scriptscriptstyle bcd}$}]
\place(320,250)[\to<125>]
\efig$$
where $\pi_1\theta_{bcd}$ is $\pi_1\cdot\pi_1\colon (X_{cd}\times Y_{cd})\cdot (X_{bc}\times Y_{bc})\to<125>X_{bd}\cdot X_{bc}$ followed by the comparison cell $X_{bd}\cdot X_{bc}\to<125>X_{bd}$, and $\pi_2\theta_{bcd}$ is defined similarly.  

\smallskip
We will say $-\times Y_{bcd}$ {\it preserves pseudo-functors} if $\theta_{bcd}$ is invertible, whenever $X$ is a pseudo-functor in $\lax(\{b,c,d\},\mathbb D)$.  In the case where $B$ is the 3-element totally-ordered set $\mathbb3=\{0,1,2\}$, we will say $-\times Y$ preserves pseudo-functors, when $-\times Y_{012}$ does.

\smallskip
Note that if $Y_{bcd}$ preserves pseudo-functors, for all $b<c<d$, then $Y$  itself is necessarily a pseudo-functor, provided that $\lax(B,\mathbb D)$ has a pseudo-functorial terminal object, e.g., $\mathbb D$ has a double terminal (in the sense of \cite{GP99}).  This is the case for the five double categories under consideration.

\theorem\label{lax} Suppose  $\mathbb D$ is a  horizontally invariant double category with a zero object such that $0\times X\cong 0$, for all $X$ in $\mathbb D$, and $\lax(\mathbb3,\mathbb D)$ has finite limits.  Then $Y\colon \mathbb3\to<125>\mathbb D$ is {\exp} in $\lax(\mathbb3,\mathbb D)$ if and only if
\description
\item{(i)} $Y_b$ is {\exp} in $\mathbb D_0$, for all $b$;
\item{(ii)} $Y_b\tov Y_c$ is {\exp} in $\mathbb D_1$, for all $b<c$; and
\item{(iii)} $-\times Y$ preserves pseudo-functors.
\enddescription
\endtheorem

\proof Suppose (i), (ii), and (iii) hold. By the remark following Lemma~\ref{subposet}, we know that $(X\times Y)_{bc}\colon X_b\times Y_b\tov X_c\times Y_c$ is the product of $X_{bc}$ and $Y_{bc}$ in $\mathbb D_1$,  for all $b<c$.  Given $Z\colon \mathbb3\to<125>\mathbb D$, consider the exponential $Z_{bc}^{Y_{bc}}$ in $\mathbb D_1$.  Then $(Z_{bc}^{Y_{bc}})_b\cong Z_b^{Y_b}$, since there are natural bijections
\begin{eqnarray*}
\mathbb D_0(X\times Y_b,Z_b) & \cong & \mathbb D_0(X\times Y_b,(Z_{bc})_0)\\
& \cong & \mathbb D_1(L_0(X\times Y_b),Z_{bc})\\
&  \cong & \mathbb D_1(L_0X\times Y_{bc},Z_{bc})\\
&  \cong & \mathbb D_1(L_0X,Z_{bc}^{Y_{bc}}) \\
&  \cong & \mathbb D_0(X,(Z_{bc}^{Y_{bc}})_0)
\end{eqnarray*}
where the third bijection follows from the isomorphism $0\times Y_c\cong 0$. Similarly, $(Z_{bc}^{Y_{bc}})_c\cong Z_c^{Y_c}$, and by horizontal invariance, we can assume these isomorphisms are equalities.

\smallskip
Consider the exponential  $Z_{bc}^{Y_{bc}}\colon Z_b^{Y_b}\tov Z_c^{Y_c}$ in $\mathbb D_1$.  Thus, $b\mapsto Z_b^{Y_b}$ becomes a lax functor $\mathbb3\to<125>\mathbb D$ via the cell on the left which corresponds to the diagram on the right by {\expty} of $Y_{02}\colon Y_0\tov Y_2$

\begin{displaymath}
\bfig
\morphism(0,600)<0,-300>[Z^{Y_0}_0`Z^{Y_1}_1;Z^{Y_{01}}_{01}]
\morphism(0,300)<0,-300>[Z^{Y_1}_1`Z^{Y_2}_2;Z^{Y_{12}}_{12}]
\square/->``->`->/<450,600>[Z^{Y_0}_0`Z^{Y_0}_0`Z^{Y_2}_2`Z^{Y_2}_2;\scriptscriptstyle id``Z^{Y_{02}}_{02}`\scriptscriptstyle id]
\place(2,200)[\hbox{\tiny$\bullet$}] 
\place(2,490)[\hbox{\tiny$\bullet$}] 
\place(450,350)[\hbox{\tiny$\bullet$}] 
\place(200,300)[\to<100>] 
\efig
\hskip .7in
\bfig
\morphism(550,600)<0,-300>[Z^{Y_0}_0\times Y_0`Z^{Y_1}_1\times Y_1;]
\morphism(550,300)<0,-300>[Z^{Y_1}_1\times Y_1`Z^{Y_2}_2\times Y_2;]
\square/->`->``->/<550,600>[Z^{Y_0}_0\times Y_0`Z^{Y_0}_0\times Y_0`Z^{Y_2}_2\times Y_2`Z^{Y_2}_2\times Y_2;\scriptscriptstyle id`\scriptscriptstyle{(Z^{Y_{12}}_{12}\cdot Z^{Y_{01}}_{01})\times Y_{02}}``\scriptscriptstyle id]
\place(2,350)[\hbox{\tiny$\bullet$}] 
\place(552,490)[\hbox{\tiny$\bullet$}] 
\place(552,200)[\hbox{\tiny$\bullet$}] 
\place(190,300)[\to<100>] 
\square(550,300)/->``->`->/<450,300>[Z^{Y_0}_0\times Y_0`Z_0`Z^{Y_1}_1\times Y_1`Z_1;{\scriptscriptstyle \varepsilon_0}`\scriptscriptstyle{Z^{Y_{01}}_{01}\times Y_{01}}`\scriptscriptstyle{Z_{01}}`]
\square(550,0)/->``->`->/<450,300>[Z^{Y_1}_1\times Y_1`Z_1`Z^{Y_2}_2\times Y_2`Z_2;{\scriptscriptstyle \varepsilon_1}`\scriptscriptstyle{Z^{Y_{12}}_{12}\times Y_{12}}`\scriptscriptstyle{Z_{12}}`{\scriptscriptstyle \varepsilon_2}]
\place(1002,490)[\hbox{\tiny$\bullet$}] 
\place(1002,200)[\hbox{\tiny$\bullet$}] 
\place(800,470)[\to<100>]  
\place(800,170)[\to<100>] 
\square(1000,0)/->``->`->/<450,600>[Z_0`Z_0`Z_2`Z_2;\scriptscriptstyle id``\scriptscriptstyle{Z_{02}}`\scriptscriptstyle id]
\place(1450,350)[\hbox{\tiny$\bullet$}] 
\place(1250,300)[\to<100>]  
\efig 
\end{displaymath}
where there is a cell in the left rectangle of the right diagram, since $-\times Y$ preserves pseudo-functors. Note that we are applying the latter to the pseudo-functor  $b\mapsto Z_b^{Y_b}$ with $Z^{Y_{12}}_{12}\cdot Z^{Y_{01}}_{01}\colon Z^{Y_0}_0\tov Z^{Y_2}_2$ and the identity comparison cell.   With this definition, it is not difficult to show that the unit and counit for $(\ )^{Y_{bc}}$ extend to ones for $(\ )^Y$, and it follows that $Y$ is {\exp} in  $\lax(\mathbb3,\mathbb D)$.

\smallskip
Conversely, suppose $Y$ is {\exp} in $\lax(\mathbb3,\mathbb D)$.  Then arguments analogous to the one proving $(Z_{bc}^{Y_{bc}})_b\cong Z_b^{Y_b}$ in the first half of the proof, show (i) and (ii) hold.  To see that  $-\times Y$ preserves pseudo-functors,  suppose $X\colon \mathbb3\to<125>\mathbb D$ is a pseudo-functor.  Then
$$\bfig
\morphism(300,800)<0,-200>[\scriptstyle0`\scriptscriptstyle{X_1};]
\place(300,715)[\cdot] 
\morphism(300,600)<0,-200>[\scriptscriptstyle{X_1}`\scriptstyle0;]
\place(300,520)[\cdot] 
\morphism(0,500)<0,-200>[\scriptscriptstyle{X_0}`\scriptscriptstyle{X_1};]
\place(0,415)[\cdot] 
\morphism(0,300)<0,-200>[\scriptscriptstyle{X_1}`\scriptscriptstyle 0;]
\place(0,220)[\cdot] 
\morphism(600,500)<0,-200>[\scriptstyle0`\scriptscriptstyle{X_1};]
\place(600,415)[\cdot] 
\morphism(600,300)<0,-200>[\scriptscriptstyle{X_1}`\scriptscriptstyle{X_2};]
\place(600,220)[\cdot] 
\Atriangle(0,300)/->`->`/<300,300>[\scriptscriptstyle{X_1}`\scriptscriptstyle{X_1}`\scriptscriptstyle{X_1};``]
\Vtriangle(0,0)/`->`->/<300,300>[\scriptscriptstyle{X_1}`\scriptscriptstyle{X_1}`X;``]
\efig$$ 
is a pushout,  by Lemma~\ref{po}. Since $-\times Y$ preserves pushouts, it follows that  
$$\bfig
\morphism(300,800)<0,-200>[\scriptstyle0`\scriptscriptstyle{X_1\times Y_1};]
\place(300,715)[\cdot] 
\morphism(300,600)<0,-200>[\scriptscriptstyle{X_1\times Y_1}`\scriptstyle0;]
\place(300,520)[\cdot] 
\morphism(0,500)<0,-200>[\scriptscriptstyle{X_0\times Y_0}`\scriptscriptstyle{X_1\times Y_1};]
\place(0,415)[\cdot] 
\morphism(0,300)<0,-200>[\scriptscriptstyle{X_1\times Y_1}`\scriptstyle0;]
\place(0,220)[\cdot] 
\morphism(600,500)<0,-200>[\scriptstyle0`\scriptscriptstyle{X_1\times Y_1};]
\place(600,415)[\cdot] 
\morphism(600,300)<0,-200>[\scriptscriptstyle{X_1\times Y_1}`\scriptscriptstyle{X_2\times Y_2};]
\place(600,220)[\cdot] 
\Atriangle(0,300)/->`->`/<300,300>[\scriptscriptstyle{X_1\times Y_1}`\scriptscriptstyle{X_1\times Y_1}`\scriptscriptstyle{X_1\times Y_1};``]
\Vtriangle(0,0)/`->`->/<300,300>[\scriptscriptstyle{X_1\times Y_1}`\scriptscriptstyle{X_1\times Y_1}`X\times Y;``]
\efig$$ 
is a pushout, as well, and the desired result follows.
\endproof

\smallskip
Note that the obstruction to extending the above proof to a general poset $B$ is that the construction $b\mapsto Z_b^{Y_b}$ may not be a lax functor since the coherence condition relative to associativity need not hold.  We will see that this is not a problem in $\cat$, or for certain double categories including the other four of interest here.

\smallskip
Using  B\'{e}nabou's equivalence $\lax(B,\cat)\simeq {\rm Cat}/B$, the fact that (i) and (ii) always hold, and $-\times Y_{bcd}$ preserves pseudo-functors, for all pseudo-functors $Y\colon B\to<125>\cat$, we get Street's characterization \cite{Street}, as a consequence of Theorem~\ref{lax} as follows.

\smallskip
\corollary A functor $Y\to<125>B$ is {\exp} in $ {\rm Cat}$ if and only if the corresponding vertical normal lax functor $B\to<125>\cat$ is a pseudo-functor.
\endcorollary

\proof Suppose $Y\to<125>B$ is {\exp} in $ {\rm Cat}$.  Then so is $\mathbb3\times_BY\to<125>\mathbb3$, for all functors $f\colon \mathbb3\to<125>B$, since pulling back along $f$ preserves {\expty}.  Thus, the corresponding normal lax functor, denoted via abuse of notation by $Y_f \colon \mathbb3\to<125>\cat$, is {\exp} in  $\lax(\mathbb3,\cat)$, and hence, satisfies (iii) of Theorem~\ref{lax}.  Since $\cat$ has a double terminal object, it follows that  that $Y_f$ is a pseudo-functor, and so $B\to<125>\cat$ is as well.

\smallskip
Conversely, suppose the corresponding vertical normal lax functor $Y\colon B \to<125> \cat$ is a pseudo-functor.  Then so is $Y_f\colon \mathbb3\to<125>\cat$, for all $f\colon \mathbb3\to<125> B$.  Since  (i) and (ii) of Theorem~\ref{lax} hold, in any case, and 
every pseudo-functor satisfies (iii), it follows that $Y_f$ is {\exp} in $\lax(\mathbb3,\cat)$.  To see that $Y$ is {\exp} in  $\lax(B,\cat)$, we need only show that the construction $b\mapsto Z_b^{Y_b}$ given in Theorem~\ref{lax} is coherent relative to associativity.  

\smallskip
For $\beta\colon b\to<100>c$ in $B$, the exponential $Z_\beta^{Y_\beta}\colon Z_b^{Y_b}\tov Z_c^{Y_c}$ can be described as follows.  Identifying $Z_b^{Y_b}$ with the usual functor category, elements of the set $Z_\beta^{Y_\beta}(\sigma_b,\sigma_c)$ correspond to cells in $\mathbb D$ of the form
$$\bfig
\square<350,350>[Y_b`Z_b`Y_c`Z_c;\sigma_b`Y_\beta`Z_\beta`\sigma_c]
\place(160,160)[\to<100>]
\place(0,190)[\hbox{\tiny$\bullet$}]
\place(350,190)[\hbox{\tiny$\bullet$}]
\efig$$
Unravelling the proof of Theorem~\ref{lax}, the comparison cells $Z_\gamma^{Y_\gamma}\cdot Z_\beta^{Y_\beta}\to<125>Z_{\gamma\beta}^{Y_{\gamma\beta}}$ are induced by the diagram
$$\bfig
\square/->`->``->/<350,700>[Y_b`Y_b`Y_d`Y_d;id`Y_{\gamma\beta}``id]
\place(0,350)[\hbox{\tiny$\bullet$}]
\place(160,350)[\to<100>]
\place(180,410)[\scriptstyle{\varphi^{\scriptscriptstyle{-1}}}]
\square(350,350)<350,350>[Y_b`Z_b`Y_c`Z_c;\sigma_b`Y_\beta`Z_\beta`]
\place(350,550)[\hbox{\tiny$\bullet$}]
\place(700,550)[\hbox{\tiny$\bullet$}]
\place(525,525)[\to<100>]
\square(350,0)<350,350>[Y_c`Z_c`Y_d`Z_d;\sigma_c`Y_\gamma`Z_ \gamma `\sigma_d]
\place(350,200)[\hbox{\tiny$\bullet$}]
\place(700,200)[\hbox{\tiny$\bullet$}]
\place(525,175)[\to<100>]
\square(700,0)/->``->`->/<350,700>[Z_b`Z_b`Z_d`Z_d;id``Z_{\gamma\beta}`id]
\place(700,350)[\hbox{\tiny$\bullet$}]
\place(890,350)[\to<100>]
\place(880,410)[\scriptstyle{\psi}]
\efig$$
where $\varphi$ and $\psi$ are the comparison cells for $Y$ and $Z$, respectively, and so coherence easily follows.  
\endproof 

\smallskip
Next, we prove a corollary that gives {\expty} results for $\pos$, $\spaces$, $\loc$, and $\rel$.   First, we recall (from  \cite{GP99}) a property of double categories which is shared by these four examples and eliminates the coherence problem in the general version of Theorem~\ref{lax}. 

\smallskip
A double category $\mathbb D$ is called {\it flat}  if its cells are determined by their domains and codomains.  In this case, there is at most one cell
$$\bfig
\square<350,350>[X_0`Y_0`X_1`Y_1;f_0`m`n`f_1]
\place(0,200)[\hbox{\tiny$\bullet$}]
\place(350,200)[\hbox{\tiny$\bullet$}]
\place(175,175)[\to<100>]
\efig$$
for all $f_0, m, n, f_1$.  

\corollary\label{laxB} Suppose $B$ is a poset, $\mathbb D$ is a flat horizontally invariant double category with a zero object such that $0\times X\cong 0$, for all $X$ in $\mathbb D$, and $\lax(B,\mathbb D)$ has finite limits.  Then $Y\colon B\to<125>\mathbb D$ is {\exp} in $\lax(B,\mathbb D)$ if and only if
\description
\item{(i)} $Y_b$ is {\exp} in $\mathbb D_0$, for all $b$;
\item{(ii)} $Y_b\tov Y_c$ is {\exp} in $\mathbb D_1$, for all $b<c$; and
\item{(iii)} $-\times Y_{bcd}$ preserves pseudo-functors, for all $b<c<d$.
\enddescription
\endcorollary

\proof
Given (i)--(iii), consider $b\mapsto Z_b^{Y_b}$, as defined in Theorem~\ref{lax}.  Since $\mathbb D$ is flat, this construction is coherent relative to associativity.  Thus, we get a normal lax functor $B\to<125>\mathbb D$, and it follows that $Y\colon B\to<125>\mathbb D$ is {\exp} in $\lax(B,\mathbb D)$.  

\smallskip
Suppose $Y\colon B\to<125>\mathbb D$ is {\exp} in $\lax(B,\mathbb D)$.  Then, by an argument analogous to the one proving $(Z_{bc}^{Y_{bc}})_b\cong Z_b^{Y_b}$ in the proof Theorem~\ref{lax}, we see that $Y_{bcd}$ is {\exp} in $\lax(\mathbb3,\mathbb D)$, for all $b<c<d$, and the desired result follows.
\endproof

\smallskip
When $B$ is a poset, B\'{e}nabou's equivalence is easily seen to restrict to $\pos$, yielding  $\lax(B,\pos)\simeq{\rm Pos}/B$.  As in  $\cat$, we know (i) and (ii) always hold (see \cite{posets}), and $-\times Y_{bcd}$ preserves pseudo-functors, for all pseudo-functors $Y\colon B\to<125>\pos$.  Applying Corollary~\ref{laxB}, we get:

\smallskip
\corollary  A morphism $Y\to<125>B$ is {\exp} in ${\rm Pos}$ if and only if the corresponding vertical normal lax functor $B\to<125>\pos$ is a pseudo-functor.
\endcorollary

\smallskip
In \cite{glue}, we showed that if $\mathbb D$ is a double category satisfying certain conditions, then  B\'{e}nabou's equivalence  generalizes to $\lax(B,\mathbb D)\simeq \mathbb D_0/\Gamma_B1$, for every finite poset $B$, where  $\Gamma_B$ is left adjoint to the constant functor  $\mathbb D_0\to<125>\lax(B,\mathbb D)$.  Examples include $\cat, \pos, \spaces$, and $\loc$.  In any case, we know pulling back preserves {\expty} in $\mathbb D_0$, and so using this equivalence, we can replace (i) and (ii) of Corollary~\ref{laxB} by the single condition that $Y_b\tov Y_c$ is {\exp} in $\mathbb D_1$, for all $b\le c$.

\smallskip
As noted above, the finiteness condition is not necessary in $\cat$ and $\pos$.  Whether it is necessary in  $\spaces$ and $\loc$ is an open question.  However, we know that $\Gamma_B1$ is the Alexandroff space on $B$ (in the sense of \cite{Alexandrov}, i.e., open sets are downward closed) in $\spaces$, and  $\Gamma_B1$ is the locale $\downarrow\!{\rm Cl}(B)$ of down-sets of $B$ in $\loc$.  Since every finite $T_0$ space is the Alexandroff space of its poset of points with the specialization order (see \cite{SS}), we get the following two corollaries:

\corollary\label{exptop} The following are equivalent for a finite $T_0$ space $B$ and a continuous map $q\colon Y\to<125>B$ with corresponding vertical  normal lax functor $n\colon B\to<125>\spaces$.
\description
\item{(a)} $q\colon Y\to<125>B$ is {\exp} in $\rm Top$.
\item{(b)} $n\colon B\to<125>\spaces$ is {\exp} in $\lax(B,\spaces)$.
\item{(c)} $Y_b\tov Y_c$ is {\exp} in $\spaces_1$, for all $b\le c$ and $-\times n_{bcd}$ preserves pseudo-functors, for all $b<c<d$.
\enddescription
\endcorollary

\smallskip
\corollary\label{exploc} The following are equivalent for a finite poset $B$ and a locale morphism $q\colon Y\to<125>\downarrow\!{\rm Cl}(B)$ with corresponding vertical  normal lax functor $n\colon B\to<125>\loc$.
\description
\item{(a)} $q\colon Y\to<125>\downarrow\!{\rm Cl}(B)$ is {\exp} in $\rm Loc$.
\item{(b)} $n\colon B\to<125>\loc$ is {\exp} in $\lax(B,\loc)$.
\item{(c)} $Y_b\tov Y_c$ is {\exp} in $\loc_1$, for all $b\le c$ and $-\times n_{bcd}$ preserves pseudo-functors, for all $b<c<d$.

\enddescription
\endcorollary

\smallskip
In the next section, we will characterize {\expty} in $\spaces_1$ and $\loc_1$, to get more complete versions of Corollaries \ref{exptop} and \ref{exploc}.  In fact, in both cases, we will also show that $-\times n_{bcd}$ preserves pseudo-functors, for all $b<c<d$, whenever $n$ is a pseudo-functor such that $n_{bc}\colon Y_b\tov Y_c$ is {\exp}, for all $b\le c$.  

\smallskip
Now,  we will use Lemma~\ref{po} to characterize pseudo-functors, when $\Gamma/B$ is an equivalence of categories.  As noted earlier,  $\pos$ satisfies this condition for any poset $B$.  Also, $\spaces$ and $\loc$ do when $B$ is finite, and $\cat$ does for any small category $B$.  Our proof below, assumes $B$ is a poset, but can be adapted to apply to $\cat$ for a general $B$.

\smallskip
Suppose $B$ is a poset and  $n\colon B\to<125>\mathbb D$.  We adopt the following abuse of notation.  Given $b<c$, let $Y_{bc}$ denote the image of $n_{bc}\colon Y_b\tov Y_c$ under $\Gamma_\mathbb2\colon\lax(\mathbb2,\mathbb D)\to<125> \mathbb D_0$, and similarly, $Y_{bcd}$ for  $b<c<d$ and $\Gamma_\mathbb3$.  If $Y\to<125>\Gamma_B1$ is the image of $n$ under $\Gamma_B$, then $Y_{bc}\cong\Gamma_B L_{bc}1\times_{\Gamma_B1}Y$ and  $Y_{bcd}\cong \Gamma_B L_{bcd}1\times_{\Gamma_B1}Y$, when  $\mathbb D_0$ is $\cat$, $\pos$, $\spaces$, and $\loc$.   In any case, one can show that  there is a commutative diagram
$$\bfig
\Atriangle(0,300)/->`->`/<300,300>[Y_c`Y_{bc}`Y_{cd};``]
\Vtriangle(0,0)/`->`->/<300,300>[Y_{bc}`Y_{cd}`Y_{bcd};``]
\efig\eqno{(3)}$$
in $\mathbb D_0$, for all $b<c<d$.   

\proposition Suppose $\mathbb D$ has $0$ and $1$, $\Gamma/\mathbb3\colon \lax(\mathbb3,\mathbb D)\to<125> \mathbb D_0/\Gamma_\mathbb31$ is equivalence of categories, and  $B$ is a poset. Then a normal lax functor $n\colon B\to<125> \mathbb D$ is a pseudo-functor if and only if  the diagram (3) is a pushout in $\mathbb D_0$, for all $b<c<d$.
\endproposition

\proof First, $n$ is a pseudo-functor if and only if $n_{bcd}\colon \mathbb 3\to<125>\mathbb D$ is, for all $b<c<d$ if and only if  the diagram (2) from Lemma~\ref{po} is a pushout, for all $b<c<d$.  Since $\Gamma/\mathbb3\colon \lax(\mathbb3,\mathbb D)\to<125> \mathbb D_0/\Gamma_\mathbb31$ is an equivalence, the latter holds if and only if (3) is a pushout, for all $b<c<d$.
\endproof

\smallskip
We conclude this section by turning our attention to $\rel$.  In this case, the functor $\Gamma/B\colon \lax(B,\mathbb \rel)\to<125> \mathbb \rel_0/\Gamma_B1$ is not an equivalence unless $B=\mathbb1$, since it is not difficult to show that $\Gamma_B1$ is a one-point set so that $ \mathbb \rel_0/\Gamma_B1\cong \rm Set$.  Thus, the theorem from \cite{glue} does not apply.  However, $\lax(B,\rel)$  is equivalent to the category ${\rm Pos}_d/B$ of posets with discrete fibers over $B$ (see \cite{lax}), and $\rel_1$ is easily seen to be cartesian closed.   Applying Corollary~\ref{laxB}, we get:

\smallskip
\corollary  Suppose $B$ is a poset.  Then $Y\to<125>B$ is {\exp} in ${\rm Pos}_d/B$ if and only if the corresponding vertical normal lax functor $B\to<125>\rel$ is a pseudo-functor.
\endcorollary

\section{Exponentiability in $\spaces_1$ and $\loc_1$} 
A space $Y$ is {\exp} in $\rm Top$ if and only if ${\cal O}(Y)$ is a continuous lattice (in the sense of Scott \cite{Scott})  if and only if $\mathbb2^Y$ exists in $\rm Top$, where $\mathbb2$ denotes the Sierpinski space $\{0,1\}$, with $\{0\}$ open but not $\{1\}$.  In this case, $\mathbb2^Y\cong{\cal O}(Y)$ with the Scott topology, which is defined as follows.

\smallskip
Recall that a subset  $H$ of a complete lattice $L$ is called {\it Scott open} if $\uparrow\!\!H=H$ and $\vee S\in H \Rightarrow  \vee F\in H$, for some finite $F\subseteq S$.   The set $\Sigma L$ of Scott open subsets is called the {\it Scott topology} on $L$.  Given $u,v\in L$, we say $u$ is {\it way below} $v$, and write $u<<v$, if $v\le\vee S\Rightarrow u\le\vee F$, for some finite $F\subseteq S$.  Then $L$ is a {\it continuous lattice} if it satisfies $v=\vee\{u\mid u<<v\}$.  A locale which is a continuous lattice is also called {\it locally compact}.  

 \smallskip 
The characterization of {\exp} spaces has appeared in many forms, but was first achieved in 1970 when Day and Kelly \cite{Day/Kelly} proved that $-\times Y$ preserves quotient maps precisely when ${\cal O}(Y)$ is a continuous lattice.  By Freyd's Special Adjoint Functor Theorem,  $-\times Y$ has a right adjoint if and only if it preserves quotient maps, for then it preserves all colimits (since coproducts are preserved in any case).  The ``technical condition" for {\expty} in ${\rm Top}/B$, proved in \cite{thesis} and referred to in the introduction, reduces to the Day/Kelly characterization when $B=1$, and has the following form when $B$ is a poset with the down-set topology. 

\smallskip
Suppose $q\colon Y\to<125>B$ is a continuous map.  Then $H\subseteq \bigsqcup_{b\in B}{\cal O}(Y_b)$ is called {\it fiberwise Scott open} provided that $H_b$ is Scott open, for all $b\in B$, and $V_c\in H_c\Rightarrow V_b\in H_b$, for all
$b< c$ and $V\in{\cal O}(Y)$.  With this topology, $\bigsqcup_{b\in B}{\cal O}(Y_b)$ becomes a space over $B$ via the projection.  Consider the (not necessarily continuous) function  $\varepsilon\colon \left(\bigsqcup_{b\in B}{\cal O}(Y_b)\right)\times_BY\to<125> \mathbb 2$ defined by $$\varepsilon(V_b,y)=\cases{0 & if $y\in V_b$ \cr 1 & if $y\not\in V_b$\cr }$$ 
Then, from \cite{thesis}, we get:

\smallskip
\lemma\label{cart} The following are equivalent for $q\colon Y\to<125>B$ in $\rm Top$.
\description
\item{(a)} $q\colon Y\to<125>B$ is {\exp} in $\rm Top$.
\item{(b)} The map $\varepsilon\colon \left(\bigsqcup_{b\in B}{\cal O}(Y_b)\right)\times_BY\to<125> \mathbb 2$ is continuous.
\item{(c)} For all $V_b\in{\cal O}(Y_b)$ and $y_b\in V_b$, there exists $H$ fiberwise Scott open such that $V_b\in H_b$ and $y_b$ is in the interior  in $Y$ of the set  \ $\bigcup_{b\in B}(\cap H_b)$.
\enddescription
\endlemma

\smallskip
We will use Lemma~\ref{cart}, in the case where $B=\mathbb2$, to show that the continuous lattice characterization of {\exp} objects in $\rm Top$ and $\rm Loc$ generalizes to $\spaces_1$ and $\loc_1$.  But, first we use this lemma  to prove that $-\times n_{bcd}$ preserves pseudo-functors, for all $b<c<d$, whenever $n$ is a pseudo-functor such that $n_{bc}\colon Y_b\tov Y_c$ is {\exp}, for all $b\le c$, thus removing the extra condition in Corollary~\ref{exptop} (and hence, in Corollary~\ref{exploc}, as well).  We begin by recalling the equivalence between $\lax(B,\spaces)$ and ${\rm Top}/B$, for a finite poset $B$.    

\smallskip
Given $m\colon B\to<125>\spaces$, write $m_{bc}\colon X_b\tov X_c$, for $b<c$.  Let $X=\bigsqcup_{b\in B}X_b$, with $U$ {\it open} if $U_b$ is open in $X_b$, for all $b$, and $U_c\subseteq m_{bc}U_b$, for all $b<c$.  A horizontal transformation $f\colon m\to<125>n$ gives rise to a continuous map $f\colon X\to<125>Y$ in the obvious way.  Note that $B$ is the space associated with the terminal object of $\lax(B,\spaces)$, and so the projection $X\to<125>B$ is continuous.  Thus, we get  $\Gamma_B\colon \lax(B,\spaces)\to<125>{\rm Top}$, which is left adjoint to the constant functor and induces an equivalence $\Gamma/B\colon \lax(B,\spaces)\to<125>{\rm Top}/B$, whose pseudo-inverse is defined as follows (see \cite{glue} for details).

\smallskip 
For $X\to<125>B$, let $m_b=X_b$, the fiber over $X$ at $b$.  Then the inclusion $i_b\colon X_b\to<125>X$ induces a locale morphism $i_b\colon {\cal O}(X_b)\to<125> {\cal O}(X)$ defined by $i_b^*=i_b^{-1}$ and ${i_b}_*(U_b)=[U_b\cup(X\setminus X_b)]^\circ$.  Given $b<c$, the composition ${\cal O}(X_b)\to<125>^{{i_b}_*} {\cal O}(X)\to<125>^{i_c^*} {\cal O}(X_c)$ is a vertical morphism in $\spaces$ which we denote by $m_{bc}\colon X_b\tov X_c$.  Thus, we get a normal lax functor $m\colon B\to<125>\spaces$, and hence a functor $\Phi_B\colon {\rm Top}/B\to<125> \lax(B,\spaces)$ which is pseudo-inverse to $\Gamma/B$.

\smallskip
Suppose $q\colon Y\to<125>B$ corresponds with $n\colon B\to<125>\spaces$.  We would like to describe the normal lax functor related to the space $\bigsqcup_{b\in B}{\cal O}(Y_b)$ with the fiberwise Scott topology over $B$.  Of course, the fiber over $b$ is ${\cal O}(Y_b)$ with the Scott topology $\Sigma{\cal O}(Y_b)$.  Using the description of open sets of $Y$ arising from the equivalence given above, one can show that $H$ is fiberwise Scott open if and only if $H_b$ is Scott open, for all $b$, and $n_{bc}U_b\in H_c\Rightarrow U_b\in H_b$, for all $b<c$, if and only if  $H_b$ is Scott open, for all $b$, and $n_{bc}^{-1}H_c\subseteq H_b$, for all $b<c$.  

\smallskip
Given $b<c$, consider $\hat n_{bc}\colon{\cal O}(Y_b)\tov {\cal O}(Y_c)$, defined by $$\hat n_{bc}H_b=\bigcup\{H_c\in \Sigma{\cal O}(Y_c)\mid n_{bc}^{-1}H_c\subseteq H_b \}$$   It is not difficult to show that $\hat n_{bc}$ preserves finite intersections, but $\hat n$ does not necessarily define a lax functor $b\mapsto  {\cal O}(Y_b)$, since $n_{bc}^{-1}H_c$ need not be Scott open, when $H_c$ is.  The latter holds precisely when $n_{bc}\colon {\cal O}(Y_b)\to<125>{\cal O}(Y_c)$ preserves directed unions.  Such a function is called {\it Scott continuous}.  

\lemma\label{nhat} If $n$ is a pseudo-functor and $n_{bc}\colon Y_b\tov Y_c$ preserves directed unions, for all $b<c$, then 
$\hat n\colon B\to<125>\spaces$ is a pseudo-functor, and $\Gamma_B\hat n=\bigsqcup_{b\in B}{\cal O}(Y_b)$, with the fiberwise Scott topology.
\endlemma

\proof First, we show that $\hat n$ is a pseudo-functor.  Since $n_{bc}\colon {\cal O}(Y_b)\to<125> {\cal O}(Y_c)$ preserves directed unions,  it is Scott continuous, and so $n_{bc}^{-1}H_c$ is Scott open, for all $H_c$ Scott open in $ {\cal O}(Y_c)$.  Then  $n_{bc}^{-1}\colon \Sigma{\cal O}(Y_b)\to<125>  \Sigma{\cal O}(Y_c)$ is left adjoint to $\hat n_{bc}$, by definition of the latter.  Since $n$ is a pseudo-functor, we know that $n_{bd}=n_{cd}n_{bc}$, and so $n^{-1}_{bd}=n^{-1}_{bc}n^{-1}_{cd}$.  Thus, it follows that $\hat n_{bd}= \hat n_{cd}\hat n_{bc}$, as desired.

\smallskip
It remains to show that $\Gamma_B\hat n=\bigsqcup_{b\in B}{\cal O}(Y_b)$ with the fiberwise Scott topology.  By definition, $H\subseteq \Gamma_B\hat n$ is open if and only if $H_b$ is open in ${\cal O}(Y_b)$, for all $b$, and $H_c\subseteq\hat n_{bc}H_b$, for all $b<c$, if and only if $H_b$ is open in ${\cal O}(Y_b)$, for all $b$, and $n^{-1}_{bc}H_c\subseteq H_b$, for all $b<c$, and the desired result follows.
\endproof

\lemma\label{expcont} If $q\colon Y\to<125>\mathbb2$ is {\exp} in $\rm Top$, then the corresponding $n\colon Y_0\tov Y_1$ preserves directed unions.
\endlemma

\proof Suppose $\{U_\alpha\}$ is directed, and consider $\bigcup nU_\alpha\subseteq n\left(\bigcup U_\alpha\right)$.  Now, ${\cal O}(Y_1)$ is a continuous lattice, since pulling back preserves {\expty}, and so given $y_1\in  n\left(\bigcup U_\alpha\right)$, there exists $V_1\in{\cal O}(Y_1)$ such that $y_1\in V_1 << n\left(\bigcup U_\alpha\right)$.  Since $q$ is {\exp}, applying Lemma~\ref{cart} with $B=\mathbb2$, there exists $H\subseteq{\cal O}(Y_0)\bigsqcup{\cal O}(Y_1)$ fiberwise Scott open such that $V_1\in H_1$ and $y_1\in W\subseteq(\cap H_0)\cup(\cap H_1)$, for some $W$ open in $Y$.  Then $H_1$ is Scott open, $V_1\in H_1$, and $V_1\subseteq n\left(\bigcup U_\alpha\right)$, and so $n\left(\bigcup U_\alpha\right)\in H_1$.  Since $H$ is fiberwise Scott open, we know $\bigcup U_\alpha\in H_0$, and so  $U_\alpha\in H_0$, for some $\alpha$.  
Also, $W_0\subseteq U_\alpha$, since $W_0\subseteq \cap H_0$.  Thus, $y_1\in W_1\subseteq nW_0\subseteq nU_\alpha\subseteq \bigcup nU_\alpha$, and it follows that $n\left(\bigcup U_\alpha\right)\subseteq \bigcup nU_\alpha$, as desired.
\endproof

\theorem\label{exptop2}  The following are equivalent for a finite $T_0$ space $B$ and a continuous map $q\colon Y\to<125>B$ with corresponding vertical  normal lax functor $n\colon B\to<125>\spaces$.
\description
\item{(a)} $q\colon Y\to<125>B$ is {\exp} in $\rm Top$.
\item{(b)} $n\colon B\to<125>\spaces$ is {\exp} in $\lax(B,\spaces)$.
\item{(c)}  $n\colon B\to<125>\spaces$ is a pseudo-functor and $Y_b\tov Y_c$ is {\exp} in $\spaces_1$, for all $b\le c$.
\enddescription
\endtheorem

\proof As before, we know (a)$\Rightarrow$(b)$\Rightarrow$(c).  To show that (c)$\Rightarrow$(a), suppose $n$ is a pseudo-functor and $Y_b\tov Y_c$ is {\exp} in $\spaces_1$, for all $b\le c$.  By Lemma~\ref{cart}, it suffices to show that  $\varepsilon\colon \left(\bigsqcup_{b\in B}{\cal O}(Y_b)\right)\times_BY\to<125> \mathbb 2$ is continuous.  By Lemmas~\ref{nhat} and \ref{expcont}, we know $\bigsqcup_{b\in B}{\cal O}(Y_b)\to<125>B$ corresponds to $\hat n\colon B\to<125>\spaces$ defined above, and so $\varepsilon$ is continuous.  Since $Y_b\tov Y_c$ is {\exp} in $\spaces_1$, for all $b\le c$, we know $Y_b$ is {\exp} in $\rm Top$, and so $\varepsilon_b\colon{\cal O}(Y_b)\times Y_b\to<125> \mathbb 2$ is continuous, for all $b$.  Thus, we have a cell
$$\bfig
\square<450,350>[{\cal O}(Y_b)\times Y_b`\mathbb2`{\cal O}(Y_c)\times Y_c`\mathbb2;\varepsilon_b`\hat n_{bc}\times n_{bc}`\idv`\varepsilon_c]
\place(0,200)[\hbox{\tiny$\bullet$}]
\place(450,200)[\hbox{\tiny$\bullet$}]
\place(250,175)[\supseteq]
\efig$$
for each $b<c$, since $n_{bc}$ is {\exp} in $\spaces_1$.  Applying $\Gamma/B$ to the associated morphism $\hat n\times n\to<125> \mathbb2$ in $\lax(B,\spaces)$, the desired result follows.
\endproof

\smallskip
To generalize continuity to vertical morphisms in $\spaces$ and $\loc$, we first note that there is a connection between the way-below relation and the Scott topology $\Sigma L$ on $L$, namely, $u<<v$ if and only if there exists $H\in \Sigma L$ such that $v\in H$ and $u\le \wedge H$ (since $u<< v$, for all $u\le \wedge H$ and $v\in H$).  It is this condition that we generalize.

\smallskip
Suppose $n\colon L_0\tov L_1$ is in $\loc$, and define $\hat n\colon \Sigma L_0\tov \Sigma L_1$ by $$\hat nH_0=\bigcup\{H_1\in\Sigma L_1 \mid n^{-1}H_1\subseteq H_0\}$$ Although $n$ is not necessarily continuous in the Scott topology, i.e., $n^{-1}H_1$ need not be Scott open when $H_1$ is, one can show that $H_1\subseteq \hat nH_0 \iff n^{-1}H_1\subseteq H_0$.  Note that this definition of $\hat n$ agrees with the one defined above for $\spaces$.

\smallskip
Given $u_1,v_1\in L_1$ and $H_0\in \Sigma L_0$, we say $u_1$ is {\it way below $v_1$ relative to $H_0$}, written $u_1<<_{H_0}v_1$, if $u_1<<v_1$ in $L_1$, $v_1\in \hat n H_0$, and $u_1\le n(\wedge H_0)$.  Then $n\colon L_0\tov L_1$ is called {\it doubly continuous} if $L_0$ is continuous and $L_1$ satisfies  $$v_1=\bigvee\{u_1\mid u_1<<_{H_0}v_1, \ {\rm for \ some}\ H_0\in \Sigma L_0\}$$  

\lemma\label{lemma2} $n\colon Y_0\tov Y_1$ is {\exp} in $\spaces_1$ if and only if $n\colon {\cal O}(Y_0)\tov {\cal O}(Y_1)$ is doubly continuous in $\loc$.
\endlemma

\proof Suppose $n\colon Y_0\tov Y_1$ corresponds to $q\colon Y\to<125>\mathbb2$ via $\spaces_1\simeq {\rm Top}/\mathbb2$.  It suffices to show that $q\colon Y\to<125>\mathbb2$ is  {\exp} in $\rm Top$ if and only if  $n\colon Y_0\tov Y_1$ is doubly continuous.

\smallskip
Suppose $q\colon Y\to<125>\mathbb2$ is  {\exp}.  Then $Y_0$ is {\exp} in $\rm Top$, since the pullback of an {\exp} map is {\exp}, and so ${\cal O}(Y_0)$ is a continuous lattice.  To see that  $n\colon {\cal O}(Y_0)\tov {\cal O}(Y_1)$ is doubly continuous, suppose $V_1\in {\cal O}(Y_1)$ and $y_1\in V_1$.  Then, by Lemma~\ref{cart}, there exists $H$ fiberwise Scott open such that $V_1\in H_1$ and $y_1\in U\subseteq(\cap H_0)\cup(\cap H_1)$, for some $U\in{\cal O}(Y)$. We claim that $U_1<<_{H_0}V_1$.  First, $U_1<<V_1$ in ${\cal O}(Y_1)$, since $U_1\subseteq\cap H_1$ and $V_1\in H_1$. Also, $U_1\subseteq n(U_0)\subseteq n(\wedge H_0)$, since $U$ is open in $Y$ and $U_0\subseteq\cap H_0$.  Finally, since $H$ is fiberwise Scott open, we know $n^{-1}H_1\subseteq H_0$, and so $V_1\in H_1\subseteq \hat nH_0$.  Thus, $U_1<<_{H_0}V_1$, as desired.

\smallskip
Conversely, suppose $n\colon {\cal O}(Y_0)\tov {\cal O}(Y_1)$ is doubly continuous.  We will show that $q\colon Y\to<125>\mathbb2$ satisfies Lemma~\ref{cart}(c).  Given $y_0\in V_0\in{\cal O}(Y_0)$, there exists $U_0<<V_0$ such that $y_0\in U_0$.  Take $H_0=\{W_0\mid U_0<<W_0\}$ and $H_1=\emptyset$.  Then $H$ is fiberwise Scott open, $V_0\in H_0$, and $y_0\in U_0\subseteq (\cap H_0)\cup(\cap H_1)$.  Given $y_1\in V_1\in{\cal O}(Y_1)$, since $n$ is  doubly continuous, there exist $H_0$ Scott open and $U_1\in{\cal O}(Y_1)$ such that $y_1\in U_1$ and $U_1<<_{H_0}V_1$.  Take $H_1=\{W_1\in \hat nH_0\mid U_1<< W_1\}$.  Then $H=H_0\sqcup H_1$ is fiberwise Scott open, since $H_1\subseteq \hat nH_0\Rightarrow n^{-1}H_1\subseteq H_0$; $V_1\in H_1$, since $V_1\in \hat nH_0$ and $U_1<<V_1$; and  $y_1\in(\wedge H_0)\cup U_1\subseteq (\cap H_0)\cup(\cap H_1)$ and $(\wedge H_0)\cup U_1$ is open, since $U_1\subseteq n(\wedge H_0)$ by definition of $U_1<<_{H_0}V_1$.  Therefore, $q\colon Y\to<125>\mathbb2$ is  {\exp} in ${\rm Top}$, as desired.
\endproof

\smallskip
In \cite{Hyland}, Hyland showed that a locale $L$ is {\exp} in $\rm Loc$ if and only if $L$ is locally compact (i.e., a continuous lattice) if and only if the exponential $S^Y$ exists in $\rm Loc$, where $S$ denotes the Sierpinski locale.   This result is constructive so it applies to internal locales in any topos, in particular, in the topos ${\rm Sh}(B)$ of set-valued sheaves on the locale $B$.  Moreover, Joyal and Tierney \cite{JT}  showed that $q\mapsto q_*\Omega_L$ sets up an equivalence between
${\rm Loc}/B$ and the category ${\rm Loc}({\rm Sh}(B))$ of internal locales in ${\rm Sh}(B)$, where $\Omega_L$ is the subobject classifier of ${\rm Sh}(L)$ and $q\colon {\rm Sh}(L)\to<125>{\rm Sh}(B)$ is the geometric morphism induced by $q\colon L\to<125>B$.  Thus,  $q\colon L\to<125>B$ is {\exp} in $\rm Loc$ if and only if $q_*\Omega_L$ is locally compact in ${\rm Loc}({\rm Sh}(B))$

\theorem\label{locthm} The following are equivalent for $n\colon L_0\tov L_1$ in $\loc$ with corresponding morphism $q\colon L\to<125>S$ in $\rm Loc$.
\description
\item{(a)} $n\colon L_0\tov L_1$ is {\exp} in $\loc_1$.
\item{(b)} $q\colon L\to<125>S$ is {\exp} in $\rm Loc$.
\item{(c)} $n\colon L_0\tov L_1$ is doubly continuous.
\item{(d)} $q_*\Omega_L$ is locally compact in ${\rm Loc}({\rm Sh}(S))$.
\enddescription
\endtheorem

\proof $(a)\Leftrightarrow(b)\Leftrightarrow(d)$ follows from $\loc_1\simeq{\rm Loc}/S\simeq{\rm Loc}({\rm Sh}(S))$.

\smallskip
\noindent $(b)\Rightarrow(c)$ Suppose $q\colon L\to<125>S$ is {\exp} in $\rm Loc$.  Then $L_0$ and $L_1$ are  {\exp} in $\rm Loc$, and so $L_0\cong{\cal O}(Y_0)$ and $L_1\cong{\cal O}(Y_1)$, for some locally compact sober spaces $Y_0$ and $Y_1$ such that $L\cong{\cal O}(Y)$, where $Y\cong \Gamma_\mathbb2 n$ and $n$ also denotes the induced vertical morphism $n\colon Y_0\tov Y_1$ in $\spaces$.  To show that $n\colon L_0\tov L_1$ is doubly continuous, by Lemma~\ref{lemma2} it suffices to show that $n\colon Y_0\tov Y_1$ is {\exp} in $\spaces_1$, or equivalently, $Y\to<125>\mathbb2$ is {\exp} in $\rm Top$.  

\smallskip
First,  we show that ${\cal O}(X)\times_S{\cal O}(Y)$ has enough points so that ${\cal O}(X)\times_S{\cal O}(Y)\cong{\cal O}(X\times_\mathbb2Y)$, for all $X\to<125>\mathbb2$.  It is easy to see that each point of ${\cal O}(X)\times_S{\cal O}(Y)$ factors through ${\cal O}(X_0)\times{\cal O}(Y_0)$ or  ${\cal O}(X_1)\times{\cal O}(Y_1)$, and the latter locales are spatial since ${\cal O}(Y_0)$ and ${\cal O}(Y_1)$ are locally compact \cite{Isbell}.  

\smallskip
Then $Y$ is {\exp} in ${\rm Top}/\mathbb2$, by Lemma~\ref{cart}, since
\begin{eqnarray*}
{\rm Top}/\mathbb2(X\times_\mathbb2 Y,\mathbb2\times\mathbb2) & \cong & {\rm Loc}/S({\cal O}(X\times_\mathbb2 Y),S\times S)\\
& \cong & {\rm Loc}/S({\cal O}(X)\times_S{\cal O}(Y),S\times S)\\
&  \cong & {\rm Loc}/S({\cal O}(X),(S\times S)^{{\cal O}(Y)}) \\
&  \cong & {\rm Top}/\mathbb2(X,{\rm pt}((S\times S)^{{\cal O}(Y)}))
\end{eqnarray*}
where $\rm pt$ is right adjoint to $\cal O$ \cite{SS}.  Thus, $n\colon Y_0\tov Y_1$ is {\exp} in $\spaces_1$, and it follows that $n\colon L_0\tov L_1$ is doubly continuous by Lemma~\ref{lemma2}.

\smallskip
$(c)\Rightarrow(d)$ Suppose $n\colon L_0\tov L_1$ is doubly continuous.  Then $L_0$ and $L_1$ are continuous lattices, and so $n\colon {\cal O}(Y_0)\tov {\cal O}(Y_1)$, for some sober spaces $Y_0$ and $Y_1$ such that $L\cong{\cal O}(Y)$, where $Y\cong \Gamma_\mathbb2 n$.  To see that $q_*\Omega_L$ is locally compact in ${\rm Loc}({\rm Sh}(S))$ it suffices to show that $q_*\Omega_Y$ is locally compact in ${\rm Loc}({\rm Sh}(\mathbb2))$, or equivalently, for all $V$ open in $Y$, $V=\vee I$, where $I$ is the ideal $I=\{U\mid U<<V\}$ in ${\rm Sh}(\mathbb2)$.  Note that $I(\{0\})=\{U_0\mid U_0<<V_0\}$ and $I(\mathbb2)=\{U\mid U<<V \ {\rm in} \ {\cal O}(Y) \ {\rm and} \ U_0<<V_0 \ {\rm in} \ {\cal O}(Y_0) \}$.  

\smallskip
Suppose $V$ is open in $Y$.  If $V_1=\emptyset$, then $V\in {\cal O}(Y_0)$, and the result follows by continuity of ${\cal O}(Y_0)$.  Otherwise, since $n\colon {\cal O}(Y_0)\tov {\cal O}(Y_1)$ is doubly continuous, for all $y_1\in V_1$, there exists $H_0$ Scott open in ${\cal O}(Y_0)$ and $U_1\in{\cal O}(Y_1)$ such that $y_1\in U_1$ and $U_1<<_{H_0}V_1$, i.e., $V_1\in \hat nH_0$ and $U_1\subseteq n(\wedge H_0)$.  Consider $U=(\wedge H_0)\cup U_1$.  It suffices to show that $U<<V$ in  ${\cal O}(Y)$ and $U_0<<V_0$ in  ${\cal O}(Y_0)$, for then $V=\vee I$, as desired.

\smallskip
To see that $U_0<<V_0$, suppose $V_0\subseteq \cup_A W_\alpha$ in  ${\cal O}(Y_0)$.  Since $V_1\subseteq n(V_0)\subseteq n(\cup_A W_\alpha)$ and $\hat nH_0$ is Scott open, we know $n(\cup_A W_\alpha)\in \hat nH_0$, and so $\cup_A W_\alpha\in H_0$ by definition of $\hat n$.  Thus, $\cup_F W_\alpha\in H_0$, for some finite $F\subseteq A$, and it follows that $U_0\subseteq \cup_F W_\alpha$, as desired.

\smallskip
To see that $U<<V$, suppose $V\subseteq \cup_A W_\alpha$ in  ${\cal O}(Y)$.  Then $V_0\subseteq \cup_A (W_\alpha)_0$, and so $U_0\subseteq \cup_{F_0} (W_\alpha)_0$, for some finite $F_0\subseteq A$, as above.  Also, $U_1\subseteq \cup_{F_1} (W_\alpha)_1$, for some finite $F_1\subseteq A$, since $U_1<<_{H_0}V_1$.  Taking $F=F_0\cup F_1$, it follows that $U\subseteq \cup_FW_\alpha$, as desired.
\endproof 

\medskip
Note that a single theorem for {\expty} in ${\rm Top}/\mathbb2$ can be obtained by combining the conditions of Lemma~\ref{cart} for $B=\mathbb2$ with Lemma~\ref{lemma2}, and adding ``$q_*\Omega_Y$ is locally compact in ${\rm Loc}({\rm Sh}(\mathbb2))$." from Theorem~\ref{locthm}.  We can also add ``$Y\colon B\to<125>\spaces$ is a pseudo-functor and $Y_b\tov Y_c$ is doubly continuous, for all $b\le c$." to the conditions in Theorem~\ref{exptop2}, and ``$Y\colon B\to<125>\loc$ is a pseudo-functor and $Y_b\tov Y_c$ is doubly continuous, for all $b\le c$." to those in Corollary~\ref{exploc}.

\begin{references*}
\bibitem[1]{Alexandrov}
P.~S. Alexandrov, \"{U}ber die Metrisation der im kleinen kompakten 
topologische R\"{a}ume, {Math. Ann.} {92} (1924), 294--301.
\bibitem[2]{Conduche}
F. Conduch\'{e}, Au sujet de l'existence d'adjoints \`{a} droite aux
foncteurs ``image r\'{e}ciproque'' dans la cat\'{e}gorie des
cat\'{e}gories, {C. R.  Acad.  Sci. Paris} {275} (1972), 
A891--894.
\bibitem[3]{Day/Kelly}
B.~J.~Day and G.~M.~Kelly, On topological quotients preserved by 
pullback or products, {Proc. Camb. Phil. Soc.} {67} 
(1970) 553--558.
\bibitem[4]{Ehres}
Charles Ehresmann, Cat\'{e}gories structur\'{e}es, Ann. Sci. \'{E}cole Norm. Sup. 80 (1963), 349--426.
\bibitem[5]{Freyd}
P.~Freyd, Abelain Categories, Harper $\&$ Row, New York, 1964.
\bibitem[6]{Giraud}
J. Giraud, M\'{e}thode de la descente, {Bull. Math. Soc. France,
Memoire} {2} (1964).
\bibitem[7]{GP99}
Marco Grandis and Robert Par\'{e}, Limits in double categories, Cahiers de Top. et G\'{e}om. Diff. Cat\'{e}g. 40 (1999),162--220.
\bibitem[8]{GP04}
Marco Grandis and Robert Par\'{e}, Adjoints for double categories, Cahiers de Top. et G\'{e}om. Diff. Cat\'{e}g. 45 (2004),193--240.
\bibitem[9]{Hyland}
J.~M.~E. Hyland, Function spaces in the category of locales,  
Springer Lecture  Notes  in  Math. {871} (1981) 264--281.
\bibitem[10]{Isbell}
J.~R. Isbell, Atomless parts of spaces, Math. Scand. 31
(1972), 5--32.
\bibitem[11]{SS}
P.~T. Johnstone, Stone Spaces, Cambridge University Press, 1982.
\bibitem[12]{JT} A. Joyal and M. Tierney, An Extension of the Galois Theory 
of Grothendieck, Amer. Math. Soc. Memoirs 309 (1984).
\bibitem[13]{CWM}
S.~MacLane,  Categories for the Working 
Mathematician, Springer-Verlag, New York,  1971.
\bibitem[14]{thesis}
S.~B. Niefield, Cartesianness, Ph.D. Thesis, Rutgers University,
1978.
\bibitem[15]{cart}
S.~B. Niefield, Cartesianness: topological spaces, uniform spaces, and
affine schemes, {J. Pure Appl. Algebra} {23} (1982), 147--167.
\bibitem[16]{posets}
S.~B. Niefield, Exponentiable morphisms: posets, spaces, locales, and 
Grothendieck toposes, {Theory Appl. Categ.} {8} (2001), 16--32.
\bibitem[17]{perfect}
S.~B. Niefield, Exponentiability of perfect maps: four approaches, {Theory Appl. Categ.} {10} (2002), 127--133.
\bibitem[18]{lax}
S.~B. Niefield, Lax presheaves and exponentiability, {Theory Appl. Categ.} {24}, (2010), 288--301
\bibitem[19]{glue}
S.~B. Niefield, The glueing construction and double categories, to appear in Proceedings of CT2010.
\bibitem[20]{Scott}
D.~S. Scott, Continuous lattices, {Springer Lecture  Notes  in  Math.} 
{274} (1972), 97--137.
\bibitem[21]{Shulman}
M. Shulman, Framed bicategories and monoidal fibrations, {Theory Appl. Categ.} 20 (2008), 650--738.
\bibitem[23]{Street} R. Street, Powerful functors, unpublished note, September 2001.
\end{references*}
\end{document}